\documentclass[12pt]{article}
\usepackage[utf8]{inputenc}
\usepackage[parfill]{parskip}
\usepackage{amsmath}
\usepackage{amsfonts}
\usepackage{amssymb}

\title{Integral representation for Euler sums of hyperharmonic numbers}
\author{Casimir Rönnlöf}
\date{}

\newtheorem{theorem}{Theorem}
\newtheorem{corollary}{Corollary}[theorem]

\begin{document}

\maketitle

\begin{abstract}
    In this short paper, we derive an integral representation for Euler sums of hyperharmonic numbers. We use results established by other authors to then show that the integral has a closed-form in terms of zeta values and Stirling numbers of the first kind. Specifically, the integral has the form of
    $$\int_0^\infty \frac{t^{m-1}\ln(1-e^{-t})}{(1-e^{-t})^r} \ dt$$
    where $m, r \in \mathbb{N}$, $m > r$ and $r\ge1$.
\end{abstract}

\section{Introduction}

 The $n$-th hyperharmonic number of order $r$, denoted by $H_{n}^{(r)}$, is defined recursively as $H_{n}^{(0)} = \frac{1}{n}$ and

$$H_{n}^{(r)}=\sum_{k=1}^{n}H_{k}^{(r-1)}\quad r,n\geq1$$ 

Using this definition, we also see that $H_n^{(1)}$ is equal to the ordinary harmonic number $H_n = \frac{1}{1}+\frac{1}{2}+\frac{1}{3}\cdots+\frac{1}{n}$. Conway and Guy \cite{conway} provided the following identity:

$$H_{n}^{(r)}=\binom{n+r-1}{r-1}(H_{n+r-1}-H_{r-1})$$

The generating function for the hyperharmonic numbers is going to come into use later in the paper:

\begin{equation}
    \sum_{n=1}^\infty H_n^{(r)}z^n=-\frac{\ln(1-z)}{(1-z)^r}
\end{equation}

In this paper we will look at how we may represent Euler sums of hyperharmonic numbers as an integral. We will use the notation $\sigma(r,m)$ to denote the Euler sum

$$\sigma(r,m) = \sum_{n=1}^\infty\frac{H_n^{(r)}}{n^m}$$

Mező and Dil \cite{mezo} established the following connection between Euler sums of hyperharmonic numbers and the Hurwitz zeta function:

\begin{equation}
    \sigma(r,m) = \sum_{n=1}^\infty H_n^{(r-1)}\zeta(m,n) \quad\quad r\ge 1, \ m > r
\end{equation}

where the Hurwitz zeta function is defined as

$$\zeta(m,n) = \sum_{k=0}^\infty \frac{1}{(m+k)^n}$$

Later, Dil and Boyadzhiev \cite{dil} found a closed form for $\sigma(r,m)$ in terms of zeta values and Stirling numbers:

\begin{equation}
    \sigma(r,m) = \frac{1}{(r-1)!}\sum_{k=1}^r \left[{r \atop k}\right]\left(\zeta_H(m-k+1)-H_{r-1}\zeta(m-k+1)+\sum_{j=1}^{r-1} \mu(m-k+1, j)\right)
\end{equation}

where $\left[{r \atop k}\right]$ is the Stirling number of first kind and 
$$\mu(r, j) = \sum_{n=1}^\infty \frac{1}{n^r(n+j)}=\sum_{k=1}^{r-1} \frac{(-1)^{k-1}}{j^k}\zeta(r+1-k)-(-1)^r \frac{H_j}{j^r}$$

\clearpage

\section{Integral representation for $\sigma(r,m)$}

The most well-known integral representation for $\zeta(m,n)$ is 

\begin{equation}
    \zeta(m,n) = \frac{1}{\Gamma(m)}\int_0^\infty \frac{t^{m-1}e^{-nt}}{1-e^{-t}} \ dt \quad\quad m>1, \ n>0
\end{equation}

where $\Gamma(m)$ is the Gamma function. This integral representation is going to come into use when proving Theorem 1, below.

\begin{theorem}
For $m, r \in \mathbb{N}$, $r\ge1$ and $m > r$, we have:
$$\sum_{n=1}^\infty\frac{H_n^{(r)}}{n^m} = -\frac{1}{\Gamma(m)}\int_0^\infty \frac{t^{m-1}\ln(1-e^{-t})}{(1-e^{-t})^r} \ dt$$
\end{theorem}

\textbf{Proof} Let us recall equation (1)

$$\sum_{n=1}^\infty\frac{H_n^{(r)}}{n^m}=\sum_{n=1}^\infty H_n^{(r-1)}\zeta(m,n)$$

Then, by plugging in the integral representation (2) for $\zeta(m,n)$, we get

\begin{equation*}
\begin{split}
\sum_{n=1}^\infty H_n^{(r-1)}\zeta(m,n) &= \sum_{n=1}^\infty H_n^{(r-1)}\frac{1}{\Gamma(m)}\int_0^\infty \frac{t^{m-1}e^{-nt}}{1-e^{-t}} \ dt \\
& = \frac{1}{\Gamma(m)} \sum_{n=1}^\infty \int_0^\infty H_n^{(r-1)} \frac{t^{m-1}e^{-nt}}{1-e^{-t}} \ dt 
\end{split}
\end{equation*}

Because the integrand is positive for all $n$ and $x$ on the intervals, according to Fubini's theorem, we can interchange the sum and integral sign

\begin{equation*}
\begin{split}
\frac{1}{\Gamma(m)} \sum_{n=1}^\infty \int_0^\infty H_n^{(r-1)} \frac{t^{m-1}e^{-nt}}{1-e^{-t}} \ dt & = \frac{1}{\Gamma(m)} \int_0^\infty \sum_{n=1}^\infty H_n^{(r-1)} \frac{t^{m-1}e^{-nt}}{1-e^{-t}} \ dt \\
& = \frac{1}{\Gamma(m)} \int_0^\infty \frac{t^{m-1}}{1-e^{-t}} \sum_{n=1}^\infty H_n^{(r-1)} e^{-nt}\ dt
\end{split}
\end{equation*}

The sum can easily be evaluated using the generating function for the hyperharmonic numbers (1)

$$\sum_{n=1}^\infty H_n^{(r-1)} e^{-nt} = -\frac{\ln(1-e^{-t})}{(1-e^{-t})^{r-1}}$$

which leaves us with

\begin{equation*}
\begin{split}
\frac{1}{\Gamma(m)} \int_0^\infty \frac{t^{m-1}}{1-e^{-t}} \sum_{n=1}^\infty H_n^{(r-1)} e^{-nt}\ dt
& =-\frac{1}{\Gamma(m)} \int_0^\infty \frac{t^{m-1}}{1-e^{-t}} \frac{\ln(1-e^{-t})}{(1-e^{-t})^{r-1}}\ dt \\
& = -\frac{1}{\Gamma(m)} \int_0^\infty \frac{t^{m-1}\ln(1-e^{-t})}{(1-e^{-t})^{r}}\ dt 
\end{split}
\end{equation*}

which completes the proof.

\begin{corollary}
Integrating by parts $m-1$ times gets us
$$\sum_{n=1}^\infty  \frac{H_n^{(r)}}{n^m}=(m-1)!\sum_{k=1}^{m}\frac{(-1)^kt^{m-k} \ \mathcal{I}_k}{(m-k)!} = \sum_{k=1}^{m} (-1)^k(m-1)_{k-1}t^{m-k} \ \mathcal{I}_k $$

where $(m-1)_{k-1}$ is the Pochammer symbol for the falling factorial and $\mathcal{I}_k$ is the $k$-th indefinite integral shown below
$$\mathcal{I}_k = \underbrace{\int\int\cdots\int}_\text{$k$ times} \frac{\ln(1-e^{-t})}{(1-e^{-t})^{r}} \ dt \cdots \ dt \ dt$$
\end{corollary}

\begin{corollary}
By recalling the generating function for the ordinary harmonic numbers
$$\sum_{n=0}^\infty H_nz^n = -\frac{\ln(1-z)}{1-z}$$

\clearpage

and then plugging that into the integral we see that $\sigma(r,m)$ could also be represented as
\begin{equation*}
\begin{split}
    \sigma(r,m) &= \frac{1}{\Gamma(m)} \int_0^\infty \frac{t^{m-1}}{(1-e^{-t})^{r-1}}\sum_{n=0}^\infty H_ne^{-nt}\ dt \\
    & = \frac{1}{\Gamma(m)}\sum_{n=0}^\infty H_n \int_0^\infty \frac{t^{m-1}e^{-nt}}{(1-e^{-t})^{r-1}} \ dt
\end{split}
\end{equation*}
\end{corollary}

\section{Closed-form for the integral}
In this section, we will look at a few examples and how the relation can help us to solve integrals in the form given before. 

\textbf{Example 1} \textit{Find the value of the integral}

$$\int_0^\infty \frac{e^tt\ln(1-e^{-t})}{e^t-1}$$

Multiplying both the denominator and numerator by $e^{-t}$ yields an integral in the form given above. Since $\Gamma(2) = 1$, we have

$$\sigma(1,2) = \sum_{n=1}^\infty \frac{H_n}{n^2} = \int_0^\infty \frac{e^tt\ln(1-e^{-t})}{1-e^t} = 2\zeta(3)$$

Notice that since $r=1$, the hyperharmonic numbers reduced to ordinary harmonic numbers.

\textbf{Example 2} \textit{Find the value of the integral}

$$\int_0^\infty \frac{t^2\ln(1-e^{-t})}{(e^{-t}-1)^2}$$

We see that the integral is in the form of $\sigma(2,3)$. Therefore, we evaluate $\sigma(2,3)$ using (3)

$$\sigma(2,3) = 2\zeta(3)+\frac{5\zeta(4)}{4}-\zeta(2)$$

Since $\Gamma(3) = 2$ we get

$$\int_0^\infty \frac{t^2\ln(1-e^{-t})}{(e^{-t}-1)^2} =  2\zeta(2)-4\zeta(3)-\frac{5\zeta(4)}{2}$$

Note that we also multiplied by $-1$ to get the answer.

\bibliographystyle{plain}
\bibliography{bibliography.bib}

\end{document}